\documentclass[acmsmall]{acmart}

\AtBeginDocument{%
  \providecommand\BibTeX{{%
    \normalfont B\kern-0.5em{\scshape i\kern-0.25em b}\kern-0.8em\TeX}}}

\setcopyright{acmcopyright}
\copyrightyear{20xx}
\acmYear{20xx}
\acmDOI{123.xx}

\acmJournal{TOMS}
\acmVolume{xxx}
\acmNumber{12}
\acmArticle{01}
\acmMonth{12}

\usepackage[ruled]{algorithm2e}

\usepackage{subcaption}

\usepackage{enumerate}
\usepackage{multirow}

\usepackage{amssymb}

\begin{document}

\title{A Parallel Direct Eigensolver for Sequences of Hermitian Eigenvalue Problems
with No Tridiagonalization}

\author{Shengguo Li}
\orcid{0000-0001-7827-6304}
\affiliation{%
  \institution{College of Computer Science, National University of Defense Technology}
  \city{Changsha}
  \country{China}
  \postcode{410073}
}
\email{nudtlsg@nudt.edu.cn}

\author{Xinzhe Wu}
\affiliation{%
  \institution{J\"{u}lich Supercomputing Centre, Forschungszentrum J\"{u}lich}
  \city{J\"{u}lich}
  \country{Germany}
  \postcode{52425}
  }
\email{xin.wu@fz-juelich.de}

\author{Jose E. Roman}
\affiliation{\institution{D. Sistemes Inform\`{a}tics i Computaci\'{o}, Universitat Polit\`{e}cnica de Val\`{e}ncia}
\country{Spain}
}
\email{jroman@dsic.upv.es}

\author{Ziyang Yuan}
\affiliation{%
  \institution{College of Science, National University of Defense Technology}
  \city{Changsha}
  \country{China}
  \postcode{410073}
}
\email{zyuan@nudt.edu.cn}

\author{Ruibo Wang}
\affiliation{%
  \institution{College of Computer Science, National University of Defense Technology}
  \city{Changsha}
  \country{China}
  \postcode{410073}
}
\email{ruibo@nudt.edu.cn}

\author{Lizhi Cheng}
\affiliation{%
  \institution{College of Science, National University of Defense Technology}
  \city{Changsha}
  \country{China}
  \postcode{410073}
}
\email{clzcheng@nudt.edu.cn}

\renewcommand{\shortauthors}{Li, et al.}

\begin{abstract}
  In this paper, a Parallel Direct Eigensolver for Sequences of Hermitian
  Eigenvalue Problems with no tridiagonalization is proposed, denoted by
  \texttt{PDESHEP}, which combines direct methods with iterative methods.
  \texttt{PDESHEP} first reduces a Hermitian matrix to its banded form, then
  applies a spectrum slicing algorithm to the banded matrix, and finally
  computes the eigenvectors of the original matrix via backtransform.
  Therefore, compared with conventional direct eigensolvers, \texttt{PDESHEP} avoids tridiagonalization,
  which consists of many memory-bounded operations.
  In this work, the iterative method in \texttt{PDESHEP} is based on the contour
  integral method implemented in FEAST.
  For the symmetric self-consistent field (SCF) eigenvalue problems,
  \texttt{PDESHEP} can be on average $1.25\times$ faster than the state-of-the-art direct solver in ELPA
  when using $4096$ processes.
  Numerical results are obtained for dense Hermitian matrices from real applications and
  large real sparse matrices from the SuiteSparse collection.
\end{abstract}



\maketitle

\section{Introduction}
\label{sec:intro}

Large-scale Hermitian eigenvalue problems arise in many scientific computing
applications (e.g. condensed matter~\cite{kohn1999nobel}, thermoacoustics~\cite{salas2015spectral}).
In the case of electronic structure calculations based on Kohn-Sham density functional
theory, a large nonlinear eigenvalue problem is iteratively solved through
the self-consistent field (SCF) procedure, in which a sequence of related linear eigenvalue
problems needs to be solved.
Depending on the discretization method and basis used, the matrices can be sparse, dense or banded.
In this work, we focus on the dense case.
Dense eigenvalue problems are usually solved by using \emph{direct methods},
and famous packages include ScaLAPACK~\cite{scalapack-book}, ELPA~\cite{elpa-library}, EigenExa~\cite{Eigenexa}, etc.
When solving sequences of related eigenvalue problems,
these eigensolvers cannot exploit any knowledge of the properties of the problem,
and thus they solve each problem of the sequence in complete isolation.

Compared to direct methods, one advantage of iterative methods is that
the approximate eigenvectors from previous SCF iterations can be used
as the initial guess to improve the convergence rate of subspace iterations (SI) for the current
SCF iteration.
Based on SI accelerated by Chebyshev polynomials,
an efficient eigensolver (ChASE) has been proposed for sequences of dense
Hermitian eigenvalue problems~\cite{winkelmann2019chase}.
Compared to Elemental's direct eigensolver~\cite{Elemental-toms}, ChASE is faster when computing
a small portion of the extremal spectrum.
ChASE becomes slower when computing a large number of eigenpairs~\cite{winkelmann2019chase}.
It is also a common problem for other iterative methods such as
the implicitly restarted Lanczos algorithm~\cite{lehoucq1998arpack} and
LOBPCG~\cite{knyazev2001toward},  since
they have to solve a projected dense eigenvalue problem as a part of the
Rayleigh-Ritz procedure.

When computing a large number of eigenpairs with iterative algorithms, the most common approach is
the spectrum slicing technique~\cite{Siam-Eigenvalue,Saad-book}, that tries to compute the eigenvalues by chunks.
Often, spectrum slicing is based on the shift-invert spectral transformation, which involves costly matrix factorizations.
Some related works include (most of them are for sparse eigenproblems):
\begin{enumerate}
  \item The shift-invert Lanczos method (SI-Lanczos).
        It was the default choice for spectrum slicing methods~\cite{aktulga2014parallel,grimes1994shifted} and
        has been used in SIPs~\cite{zhang2007sips}, SLEPc~\cite{campos2012strategies}, and SIESTA-SIPs~\cite{kecceli2016shift,kecceli2018siesta}.

  \item Polynomial filtering based method. It uses matrix polynomial filters to amplify the spectral
        components associated with the interested eigenvalues~\cite{banerjee2016chebyshev,li2016thick}, avoiding matrix factorizations,
        and packages include EVSL (Eigenvalue Slicing Library)~\cite{li2019eigenvalues} and ChASE~\cite{winkelmann2019chase}.

  \item Contour integral based method, where a matrix factorization is required at each integration point.
  The well-known packages include FEAST~\cite{kestyn2016pfeast,polizzi2009density,peter2014feast} and z-Pares~\cite{SS-Projection,sakurai2007cirr}.

  \item The shift-invert subspace iteration method. A new method (SISLICE) is recently proposed
        in~\cite{williams2020shift} which combines shift-invert subspace iteration with spectrum
        slicing.
\end{enumerate}

Iterative methods are rarely involved to solve dense Hermitian eigenvalue problems
since both dense matrix vector multiplication and solving a dense linear system are considerably expensive.
In this paper, iterative methods are combined with direct methods in a novel way.
The central idea is that iterative methods are applied on a banded matrix with small bandwidth (usually less than $64$)
instead of the original matrix.
Then, iterative methods can be used to compute partial or full eigendecomposition of a dense Hermitian or symmetric matrix.
It is a new framework for solving dense Hermitian eigenvalue problem, and
 will be named \texttt{PDESHEP} in this paper, short for \emph{a Parallel Direct Eigensolver for
Sequences of Hermitian Eigenvalue Problems with no tridiagonalization}.

\begin{figure}[t]
\centering
\includegraphics[width=4.5in,height=1.3in]{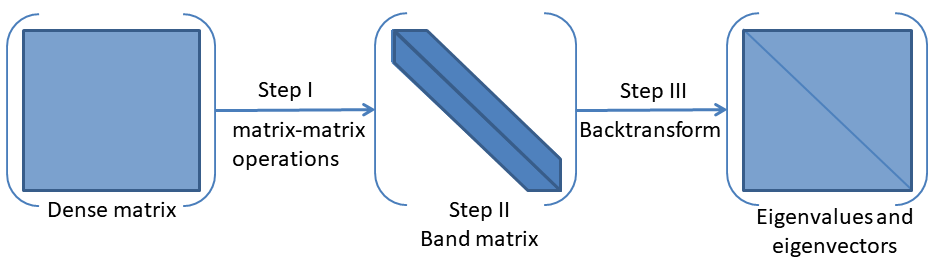}
\Description{The procedure of our algorithm}
\caption{The main steps of direct methods when using two-step symmetric banded reduction.}
\label{fig:shift-invert_subspace}
\end{figure}

The classical direct methods consist of \emph{three steps}, see Fig.~\ref{fig:shift-invert_subspace} for
the main procedure.
First, a dense matrix is reduced to its tridiagonal form via orthogonal transformations.
One efficient reduction algorithm for tridiagonalization is the two-step approach~\cite{Bischof-twostep},
which first reduces a dense matrix to band form and then further reduces it to a tridiagonal matrix.
The symmetric band reduction (SBR) from band to tridiagonal fully consists of
memory-bound operations, and its scalability is limited.
Then, the eigendecomposition of tridiagonal matrix is computed by classical algorithms and
finally the eigenvectors are computed via backtransform.
The details are introduced in the next section.
The starting point of this work is trying to get around SBR and propose a more efficient algorithm
to solve a symmetric banded eigenvalue problem.
From the results of Example 1 in section~\ref{sec:num}, we see that step (II) takes about \emph{half} of total time
when using many processes.
It is better to replace step (II) by a more efficient and faster solver.

Many eigenvalue algorithms have been proposed for symmetric banded matrices.
The most commonly used one is to reduce the banded matrix to tridiagonal form and then apply the DC, QR or
MRRR~\cite{MRRR_TOMS} algorithms to it, just as done in ELPA~\cite{Elpa}.
Some works suggest to use the banded DC algorithm~\cite{Arbenz-bdc,Gansterer-toms,Liao-camwa} when the bandwidth is small, which requires more floating point operations.
Just as shown in~\cite{thpc_pbsdc}, the banded DC algorithm is slower than ELPA when the bandwidth is large.
Another completely different approach is to use iterative methods.
In this work, we apply the spectrum slicing methods to the symmetric banded matrix instead
of the original full (dense) matrix.
Since the semibandwidth of the intermediate banded matrix is very small, $n_{bw}=64$,
it is very cheap to solve a banded linear system.
Our method is similar to~\cite{williams2020shift,zhang2007sips}, but
the matrix in our problem is banded, which makes the linear systems easy to solve.
Similarly to previous methods (SIESTA-SIPs, SISLICE, ChASE), our
algorithm is designed for SCF eigenvalue problems, and prior knowledge of eigenvalues and
eigenvectors are used for spectrum partition and as initial vectors, respectively.
The main advantage of \texttt{PDESHEP} is that it is suitable for many eigenpairs of dense matrices,
not only for sparse matrices.

\begin{figure}[b]
\centering
\includegraphics[width=4.2in,height=0.5in]{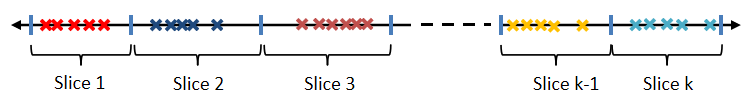}
\caption{Partitioning of the spectrum with local clusters.}
\label{fig:slice_partition}
\end{figure}

For spectrum slicing algorithms, the spectrum partition is very important for
efficiency and robustness. It relates to load balance and the orthogonality of
computed eigenvectors.
The ideal case is that every slice has nearly the same number of eigenvalues and
the gaps between eigenvalues in different slices are large, as shown in
Fig.~\ref{fig:slice_partition}.
Since we are dealing with self-consistent field (SCF) eigenvalue problems,
we assume the prior knowledge of eigenvalue distributions is known.
Otherwise, the spectral density estimation (also known as density of states or DOS)~\cite{lin2016approximating} and
Sylvester's inertia theorem~\cite{sylvester1852xix} can be used to estimate or compute
the distribution of eigenvalues of a banded matrix.
In this work, we use a simplified \emph{k}-means algorithm to partition the spectrum,
that exploits the one dimensional structure of eigenvalues which are stored increasingly.
Therefore, the simplified version requires much fewer operations than the general \emph{k}-means algorithms.
We first partition eigenvalues equally into $K$ groups and use their centers as the initial conditions
of our simplified \emph{k}-means algorithm, where $K$ is the number of partitioned slices.
We find that this initialization method for \emph{k}-means algorithms usually leads to a more `balanced' partition,
i.e, the variance of the number of eigenvalues in all slices is small.
The simplified \emph{k}-means algorithm is introduced in section~\ref{sec:spectrum_partition}.

To combine direct methods with iterative methods, we need some efficient data redistribution routines since
the direct methods such as in ScaLAPACK and ELPA store matrices in 2D block cyclic data distribution (BCDD) while
the iterative methods such as in FEAST store matrices in 1D block data distribution (BDD).
We need two efficient data redistribution algorithms. One is used to redistribute a banded
matrix in BCDD to 1D compact form as in LAPACK.
The other one is used to redistribute a dense matrix stored in \emph{irregular} 1D BDD to 2D BCDD with given
block size $n_b$.
By \emph{irregular} we mean different processes may have different number of block rows or columns.
None of BLACS routines works for the irregular case.
We used the data redistribution algorithm proposed in~\cite{tpds_dataredist}.
For the case that BLACS fits, the new data redistribution algorithm can be more than $10\times$ faster
than the BLACS routine \texttt{PXGEMR2D}.

\texttt{PDESHEP} can be seen as a new framework of direct methods, and it can be combined with any spectrum slicing methods such as
FEAST, SIPs or EVSL.
In our current version, it is implemented by combining ELPA~\cite{elpa-library} with FEAST~\cite{kestyn2016pfeast}.
Later, we will try to include SLEPc~\cite{campos2012strategies}, EVSL~\cite{li2019eigenvalues} and others.
In summary, the main contributions of this work are:
\begin{enumerate}
  \item We propose a new framework for solving the (dense) Hermitian or symmetric eigenvalue problems,
        which combines direct and iterative methods, and is suitable for computing many eigenpairs.


    \item We propose a simplified 1D \emph{k}-means algorithm for partitioning the spectrum, which only requires
    $O(n+K)$ floating operations, faster than the general \emph{k}-means algorithm, where $n$ is the number of
    points and $K$ is the number of clusters.

  \item We perform numerous experiments to test \texttt{PDESHEP} and use some dense Hermitian matrices from real applications and
        real large sparse matrices from the SuiteSparse collection~\cite{davis2011university}.

\end{enumerate}

The remainder of this paper is organized as follows.
In Section~\ref{sec:prelim}, we briefly review the main steps of direct methods for
solving dense symmetric or Hermitian eigenvalue problems and
describe the iterative methods based on contour integration, especially
the one implemented in FEAST.
In section~\ref{sec:proposed}, we present our new algorithm (\texttt{PDESHEP})  for
sequences of symmetric or Hermitian eigenvalue problems, and some related techniques such as
spectrum partition, data redistribution and validation are included.
We present numerical results in section~\ref{sec:num} and future works in section~\ref{sec:conclusion}.

\section{Preliminaries}
\label{sec:prelim}

In this section, we give a brief introduction to the direct methods and
contour integral based iterative methods for eigenvalue problems.
As mentioned above, classical direct eigensolvers can not exploit any knowledge of the properties of
related eigenvalue problems.
Thus, the algorithms in this section are only for a single eigenproblem and
how they are combined and used for solving sequences of related
eigenproblems will be introduced in section~\ref{sec:proposed}.

\subsection{Direct Methods for (Dense) Eigenvalue Problems}

In this section, we introduce the main steps of direct methods
for solving a standard eigenvalue problem,
\begin{equation}
\label{eq:evp}
AX = X\Lambda,
\end{equation}
where $A$ is an $n\times n$ Hermitian or real symmetric matrix,
$X$ is the $n\times nev$ eigenvector matrix and the
diagonal matrix $\Lambda$ contains eigenvalues $\lambda_i$, $i=1,\ldots,nev$.
Throughout this paper we will focus on the Hermitian case, but
our algorithm is applicable to real symmetric matrices.
Furthermore, the generalized eigenvalue problem $AX=\lambda BX$, where
$A$ is Hermitian and $B$ is Hermitian and positive definite, can be reduced
to the standard case via Cholesky factorization,
\begin{equation}
A X=B X \Lambda  \Longrightarrow \widetilde{A} \widetilde{X} = \widetilde{X} \Lambda,
\end{equation}
where $B= LL^H$, $\widetilde{A} = L^{-1} A L^{-H}$ and $\widetilde{X} = L^H X.$

Only the standard eigenvalue problem~\eqref{eq:evp} is considered in this paper.
We mainly follow the description in~\cite{elpa-library} and
introduce the main stages of a direct method,
which consists of \emph{three} steps.
The input dense matrix is first reduced to tridiagonal form
by using \emph{one-step} or \emph{two-step} reduction methods~\cite{Golub-book2,Bischof-twostep}.
The one-step method transforms a dense matrix to tridiagonal form
through sequences of Householder transformations directly.
In contrast, the two-step method~\cite{Bischof-twostep} first reduces a dense matrix to band form and
then further reduces it to tridiagonal form, which is shown to be more effective~\cite{Elpa}.
We only introduce here the \emph{two-step} reduction method.
\begin{enumerate}[Step (I):]

  \item Reduce $A$ to band form through a sequence of unitary
        transformations,
        \begin{equation}
        \label{eq:densetoband}
            D = U {A} U^H,
        \end{equation}
        where $U$ is a unitary matrix.
        This step relies on efficient matrix-matrix operations~\cite{Bischof-twostep},
        and can be computed efficiently on high performance computers.

  \item Compute the eigendecomposition of the banded matrix $D$,
        \begin{equation}
        \label{eq:bandsolver}
            D \widehat{X} = \widehat{X} \Lambda.
        \end{equation}
        This is achieved by the following three steps.

        \begin{enumerate}[\hspace{0.7cm}]
          \item[(IIa)] Reduce the banded matrix $D$ to tridiagonal form,
                       \begin{equation}
                        T = VD V^H,
                       \end{equation}
                       where $V$ is a unitary matrix.
                       This is done via the \emph{bulge chasing} procedure~\cite{Parlett-book,acm807},
                       and the operations are memory-bounded.

          \item[(IIb)] Solve the tridiagonal eigenvalue problem:
                     \begin{equation}
                       T \widehat{Y} = \widehat{Y} \Lambda.
                      \end{equation}
                      It can be solved by QR~\cite{francis1962qr},                       divide-and-conquer (DC)~\cite{Cuppen81,Gu-eigenvalue,Tisseur-DC},
                       Multiple Relatively Robust Representation (MRRR)~\cite{Dhillon-thesis,willems2013framework,Petschow-sisc}.
                      ELPA implements an efficient DC algorithm~\cite{Elpa}.

          \item[(IIc)] Backtransform to obtain the eigenvectors of the banded matrix,
                       \begin{equation}
                        \widehat{X} = V^H \widehat{Y}.
                       \end{equation}
        \end{enumerate}

  \item Backtransform to get the eigenvectors ${X}$ of matrix $A$,
        \begin{equation}
        \label{eq:backtransform}
            {X} = U^H \widehat{X}.
        \end{equation}

\end{enumerate}

The steps (I)-(III) are the classical approach for solving dense Hermitian
eigenvalue problems.
Some other approaches have been proposed to compute the eigendecomposition of the intermediate
banded matrix.
When the bandwidth is small, it may be more efficient
to use the banded DC (BDC) algorithm~\cite{Haidar-sisc,Liao-camwa,thpc_pbsdc},
but it is very slow when the bandwidth is large~\cite{thpc_pbsdc}.
The package EigenExa~\cite{Eigenexa} also implements the BDC algorithm working on a
pentadiagonal matrix (band matrix with semibandwidth $2$).

Figure~\ref{fig:shift-invert_subspace} shows the main steps of direct methods.
Step (II) is traditionally based on tridiagonalization via symmetric bandwidth reduction (SBR)~\cite{acm807} and
solved by DC or MRRR algorithms.
The starting point of this work is that step (II) has poor scalability and it takes a large portion of
total time.
See Table~\ref{tab:elpa_step4} in section~\ref{sec:num} for more details.
In this work, we instead replace step (II) by a spectrum slicing method, and
steps (I) and (III) are the same as traditional direct methods.
The advantage is that we get around SBR, spectrum slicing methods are naturally parallelizable and hence \texttt{PDESHEP} is more efficient than classical
direct methods when using many processes.


\subsection{The FEAST Algorithm}

The FEAST algorithm utilizes spectral projection and subspace iteration to obtain
selected interior eigenpairs.
Note that we still use the standard eigenvalue problem~\eqref{eq:evp} to introduce FEAST,
which works for generalized eigenvalue problems too.
A Rayleigh-Ritz procedure is used to project matrix $A$ onto a reduced search subspace to
form matrix
\begin{equation}
A_q = Q^H AQ.
\end{equation}
Approximate eigenvalues $\widetilde{\Lambda}$ and eigenvectors $\widetilde{X}$ of the original
matrix can be obtained from the much smaller eigenvalue problem,
\begin{equation}
A_q W_q = W_q \Lambda_q,
\end{equation}
as
\begin{equation}
\widetilde{X} = Q W_q, \quad \widetilde{\Lambda} = \Lambda_q.
\end{equation}

In theory, the ideal filter for the Hermitian problem would act as a projection operator $X_m X_m^H$ onto
the subspace spanned by the eigenvector basis, which can be expressed via the Cauchy integral formula
\begin{equation}
\label{eq:integral}
X_m X_m^H = \oint_\Gamma  (zI-A)^{-1} dz,
\end{equation}
where the eigenvalues associated with the orthonormal eigenvectors $X_m$ are
located within a interval delimited by the closed curve $\Gamma$.
In practice, the spectral projector is approximated by a quadrature rule using $n_e$
integration nodes and weights $\{(z_j, w_j)\}_{j=1,\ldots,n_e}$,
\begin{equation}
Q_{m_0} = \sum_{j=1}^{n_e} w_j (z_jI-A)^{-1}\widetilde{X}_{m_0} \triangleq r_m(A) \widetilde{X}_{m_0},
\end{equation}
where the search subspace is of size $m_0 \ge m$, and
$r_m(x) = \sum_{j=1}^{n_e} \frac{w_j}{z_j -x}$ is a rational function of order $(n_e-1,n_e)$.
The computation of $Q_{m_0}$ amounts to solving a set of $n_e$ independent complex
shifted linear systems with multiple right-hand sides,
\begin{equation}
\label{eq:feast_L2}
(z_j I-A)Q_{m_0}^{(j)} =  \widetilde{X}_{m_0}, \qquad \text{ for } j=1,\ldots, n_e,
\end{equation}
and $Q_{m_0} = \sum_{j=1}^{n_e} w_j Q_{m_0}^{(j)}.$
The matrix $Q_{m_0}$ is then used as the Rayleigh-Ritz projector to form the
reduced matrix $A_q$.
The general outline of FEAST is shown in Algorithm~\ref{alg:FEAST} for
computing $m$ eigenpairs in a given interval.
The initial vectors $\widetilde{X}_{m_0}$ are usually a set of random vectors, and
good initial vectors can accelerate the convergence~\cite{peter2014feast}.
It is mathematically equivalent to subspace iteration applied with a rational
matrix function $r_m(A)$, see~\cite{peter2014feast}.
The rational function can be constructed from numerical approximation of the integral~\eqref{eq:integral}
by using Gauss quadrature~\cite{polizzi2009density} or Trapezoid rule~\cite{SS-Projection}.
An improved rational approximation based on the work of Zolotarev is proposed
in~\cite{guttel2015zolotarev}, which improves both convergence robustness and load balancing
when FEAST runs on multiple search intervals in parallel.


\begin{algorithm}[thb]
\DontPrintSemicolon
\KwIn{The input matrix $A$, initial vectors $\widetilde{X}_{m_0}$,
  the residual tolerance $\epsilon$, integration nodes and weights $\{(z_j, w_j)\}_{j=1,\ldots,n_e}$.}
\KwOut{The approximated eigenvalues and eigenvectors: $\widetilde{X}_m$, $\widetilde{\Lambda}_m$.}

\While{$( \|A\widetilde{X}_m -\widetilde{X}_m \Lambda_m \| \ge \epsilon)$
    }{
      $Q_{m_0} = 0$\;
       Compute $Q_{m_0} = r_m(A)\widetilde{X}_{m_0}$\;
       construct $A_q = Q_{m_0}^H A Q_{m_0}$, and solve $A_q W_q = W_q \Lambda_q$\;
       compute $\widetilde{X}_{m_0} = Q_{m_0}W_q, \widetilde{\Lambda} = \Lambda_q$
        }

\caption{The FEAST Algorithm}
\label{alg:FEAST}
\end{algorithm}

\begin{figure}[hbt]
\centering
\includegraphics[width=3.5in,height=1.4in]{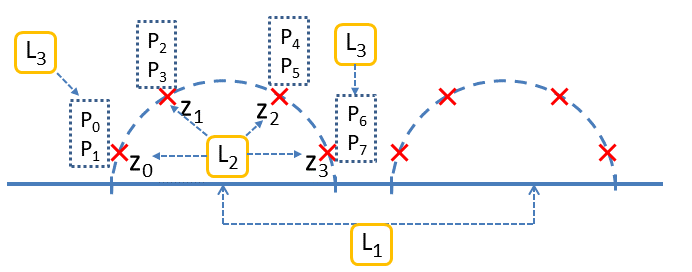}
\Description{The procedure of our algorithm}
\caption{The three levels of parallelism inherent to FEAST.}
\label{fig:pfeast_3level}
\end{figure}

There are three levels of parallelism inherent to the FEAST algorithm, which are
denoted as \textbf{L1}, \textbf{L2} and \textbf{L3}.
Figure~\ref{fig:pfeast_3level} shows the three-level parallelism of FEAST.
The first level \textbf{L1} parallelism divides the search interval into several
slices and the eigenvalues in each slice are computed by a subgroup of processes.
The second level \textbf{L2} parallelism means that each linear system in equation~\eqref{eq:feast_L2}
can be solved independently by different groups of processes.
The third level \textbf{L3} parallelism means that each linear system in
equation~\eqref{eq:feast_L2} can be solved in parallel by using multithreads or multiple processes.
The total number of MPI processes that can be used equals the product of
number of slices in \textbf{L1}, number of shifted linear systems solved simultaneously in \textbf{L2} and number of MPI
used to solve each linear system in \textbf{L3}, i.e.,
\begin{equation}
\label{eq:mpi_numbers}
\Xi(slices) \times \Xi(linear\_systems\_solved\_simultaneously) \times \Xi(MPI\_for\_one\_linear\_system),
\end{equation}
where $\Xi(x)$ means the number of entries in set $x$.

The parallelism of \textbf{L2} is limited by the number of integration nodes.
FEAST uses eight quadrature points per interval by default,
and each quadrature point corresponds to a shifted linear system.
However, the numerical results in~\cite{aktulga2014parallel} show that
using $16$ quadrature points per interval usually gives the best overall performance
from all the tested cases (8, 10, 12, 16, 20 and 24 quadrature points).
Therefore, we limit the parallelism of \textbf{L2} to $16$ in our implementation.
Different MPI process groups of \textbf{L2} solve different shifted linear system simultaneously.
For example, we can use 16 process groups to solve 16 shifted linear systems simultaneously if $n_e=16$.

The \textbf{L3} level of parallelism is not used in our experiments, since
FEAST-4.0\footnote{We downloaded from \url{http://www.ecs.umass.edu/~polizzi/feast/download.htm}.}
 does not fully support
\textbf{L3} parallelism for banded eigenvalue problems.
Each linear system in \textbf{L3} is solved by using one MPI process in our numerical
experiments.
Because we currently cannot use \textbf{L3} parallelism and
the number of MPI processes used in \textbf{L2} is limited to $16$,
the spectrum must be partitioned into many slices when using many processes.
The spectrum is partitioned into up to $256$ slices in our experiments when using
$4096$ processes,
see Example 3 in section~\ref{sec:num}.
Future work includes adding support for a distributed-memory parallel solver for banded linear systems.

\section{Proposed algorithm}
\label{sec:proposed}

In this section, we introduce our new algorithm \texttt{PDESHEP},
which can be seen as a hybrid algorithm, that combines the techniques of direct
and iterative methods.

To reduce computational complexity, the direct methods first reduce a dense matrix $A$
to a canonical form such as tridiagonal, bidiagonal or Hessenberg, and then apply QR,
DC or MRRR to compute the eigenvalue or singular value decomposition.
\texttt{PDESHEP} also follows this approach.
It first reduces a dense Hermitian matrix to its banded form
via a sequence of unitary transformations.
In contrast to classical direct algorithms, \texttt{PDESHEP} uses spectrum slicing methods to compute
the eigendecomposition of the banded matrix.
After obtaining the eigenpairs of the banded matrix, the eigenvectors of $A$ are computed via
backtransforms.
The advantage is that it gets around the SBR process, which
consists of memory-bound operations and has poor scalability just as shown in Example 1 in section~\ref{sec:num}.

\begin{figure}[hbt]
\centering
\includegraphics[width=4.0in,height=2.3in]{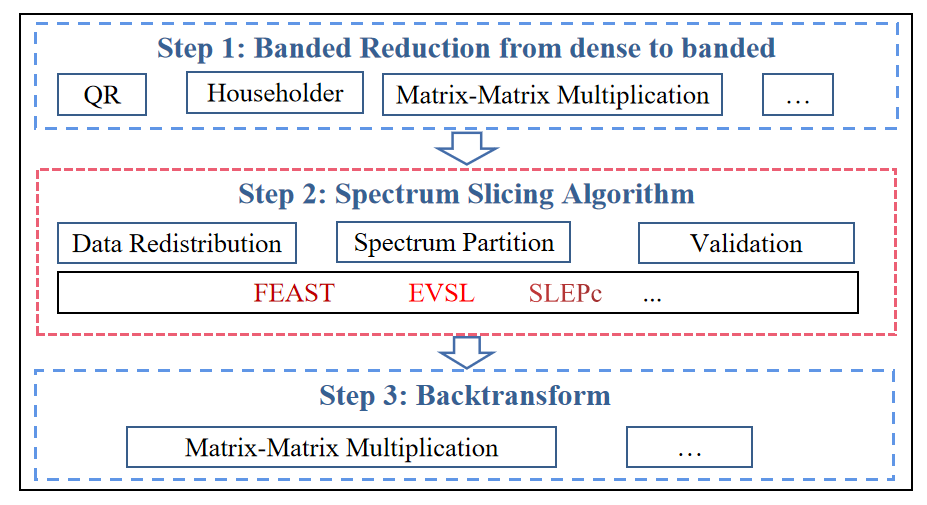}
\caption{The structure of PDESHEP.}
\label{fig:NewFramework}
\end{figure}

The procedure of \texttt{PDESHEP} is shown in Algorithm~\ref{alg:PDESHEP} and see also Figure~\ref{fig:NewFramework}.
The structure of \texttt{PDESHEP} looks like a \emph{sandwich}:
the two outer layers remain unchanged as in the classical direct methods, but the inner kernel
is replaced by an iterative method.
\texttt{PDESHEP} is more like a framework, and we do not specify which spectrum slicing algorithm to use.
Step 2 of \texttt{PDESHEP} can use any of the four types of spectrum slicing algorithms mentioned
in section~\ref{sec:intro}.
In Algorithm~\ref{alg:PDESHEP}, the input initial vectors
$\widehat{X}_{m_i}^{(\ell-1)}$ are constructed from the eigenvectors of the banded
matrix reduced from previous (dense) Hermitian matrix $A^{(\ell-1)}$,
which may be interlaced with some random vectors when more initial vectors 
are required.
Note that $\widehat{X}_{m_i}^{(\ell-1)}$ are the eigenvectors of the intermediate banded matrix instead of the original dense
Hermitian matrix.
The classical dense numerical linear algebra packages store matrices in 2D BCDD format
such as ScaLAPACK and ELPA.
The standard iterative methods usually store matrices in 1D BDD  format such as
FEAST~\cite{kestyn2016pfeast} and SLEPc~\cite{slepc}.
We need one more data redistribution if using $X$ as the initial eigenvectors.

\begin{algorithm}[thb]
\DontPrintSemicolon
\KwIn{Hermitian matrix $A^{(\ell)}$, number of eigenpairs to be computed $nev$, initial vectors $\widehat{X}_{m_i}^{(\ell-1)}$ for the intermediate banded matrix.}
\KwOut{The computed eigenpairs $(X, \Lambda)$, $X \in C^{n \times nev}$ and the
       intermediate eigenvectors $\widehat{X}_{m_i}^{(\ell)}$.}

\textbf{Step (I): Banded reduction} \\
 \begin{itemize}
   \item Apply block Householder transformations to $A^{(\ell)}$ and reduce it to a banded matrix,
   \[
   D^{(\ell)} = U {A^{(\ell)}} U^H, \quad \text{ see equation~\eqref{eq:densetoband}.}
   \]
 \end{itemize}

 \textbf{Step (II): Spectrum slicing algorithm} \\
   \begin{itemize}
     \item Redistribute the banded matrix to each process and store it locally in compact form.\;

     \item Apply a spectrum slicing algorithm to matrix $D^{(\ell)}$ by using $\widehat{X}_{m_i}^{(\ell-1)}$ as
           the initial vectors, and
           compute $nev$ eigenpairs in parallel,
           \[
           D^{(\ell)} \widehat{X} = \widehat{X} \Lambda, \quad \text{ see equation~\eqref{eq:bandsolver}.}
           \]

     \item Redistribute the computed eigenvectors $\widehat{X}$ to 2D BCDD form $\widetilde{X}$.\;
   \end{itemize}

\textbf{Step (III): Backtransform}
 \begin{itemize}
   \item Compute the eigenvectors $X$ of matrix $A^{(\ell)}$ via backtransform,
   \[
      {X} = U^H \widetilde{X}, \quad \text{ see equation~\eqref{eq:backtransform}.}
   \]
 \end{itemize}

\caption{\texttt{PDESHEP} for sequences of Hermitian eigenvalue problems}
\label{alg:PDESHEP}
\end{algorithm}

In this work, the banded eigenvalue problem is solved by the contour integral based method implemented in
FEAST~\cite{kestyn2016pfeast}.
We currently focus on the Hermitian standard eigenvalue problem.
For the generalized eigenvalue problem, we assume it has been reduced to the standard form.
FEAST can be seen as a subspace iteration eigensolver accelerated
by approximate spectral projection~\cite{peter2014feast}, and
good starting vectors $\widehat{X}_{m}$ can speed up the convergence rate.
Since the input matrices are correlated in our problem,
we use the eigenvectors of the previous banded intermediate matrix $D^{(\ell-1)}$ as the
starting vectors for the current banded matrix $D^{(\ell)}$.
The intuition is that as the SCF converges, matrices $\{A^{(\ell)}\}_{\ell\ge \eta}$ change very little,
and so do the intermediate banded matrices $\{D^{(\ell)}\}_{\ell \ge \eta}$
after $\eta$ SCF iterations, where $\eta$ is an integer number.
Numerical results verify this fact and show that using intermediate eigenvectors computed in the previous SCF iteration as starting vectors
can significantly improve the convergence, which is usually more than $1.50$ times faster, and the
results are shown in Table~\ref{tab:wo_vectors} in section~\ref{sec:num}.

\begin{figure}[ptbh]
\centering
\subcaptionbox{ The procedure of redistributing a banded matrix to 1D compact form.\label{fig:2Dto1D}}{
\includegraphics[width=2.5in,height=1.1in]{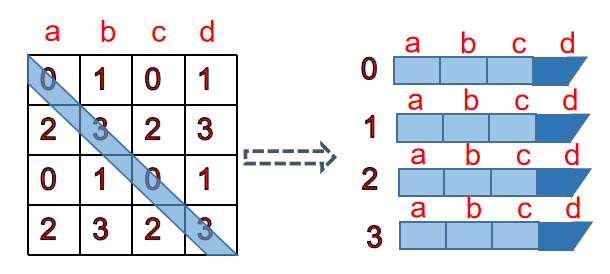}
}
\subcaptionbox{Data redistribution from 1D irregular block to 2D BCDD form.\label{fig:1Dto2D}}{
\includegraphics[width=2.5in,height=1.1in]{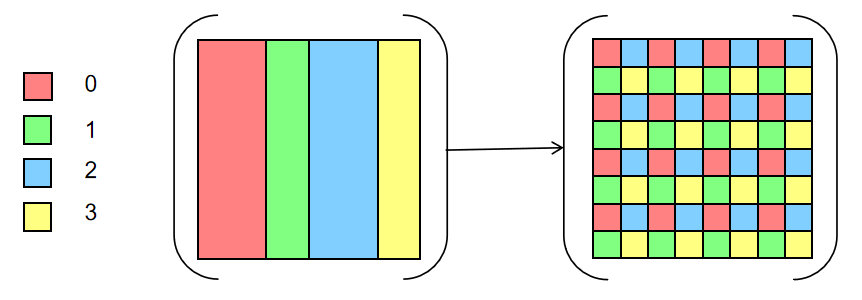}
}
\caption{The data redistribution algorithms that glue direct and iterative methods.}
\label{fig:data_redistribution}
\end{figure}

The dense packages such as ScaLAPACK and ELPA  store the matrices in 2D BCDD form.
In contrast, iterative methods, which are usually developed for sparse matrices, use
1D block data distribution to store matrices.
Since the semibandwidth of the intermediate banded matrix $D$ is small,
usually $n_{bw}=64$, we store the full $D$ locally on each process.
Therefore, we first need a data redistribution to obtain $D$ by collecting data from all
processes, which is originally stored in 2D BCDD.
The process is shown in Fig.~\eqref{fig:2Dto1D}, and
then $D$ is stored compactly as in LAPACK.
Later, we need to redistribute the computed eigenvectors of the banded matrix
from 1D BDD to 2D BCDD after applying the spectrum slicing algorithm.
This process is shown in Fig.~\eqref{fig:1Dto2D}, and
the details of data redistribution are described in~\cite{tpds_dataredist}.

\subsection{Spectrum Partition}
\label{sec:spectrum_partition}

In spectrum slicing algorithms, the spectrum is partitioned into several slices, and
the eigenvalues in each slice are computed by a group of processes.
It corresponds to the first level \textbf{L1} parallelism.
One classical way of spectrum partition is based on Sylvester's law
of inertia, which computes the $LDL^H$ factorization of $A-s_i I$~\cite{Golub-book2}.
A more modern approach is to estimate the \emph{density of states} (DOS)~\cite{lin2016approximating,xi2018fast}
of $A$, which relies on matrix-vector multiplications.
The DOS has been used to partition the spectrum for the first SCF iteration in~\cite{williams2020shift}
and EVSL~\cite{li2019eigenvalues}.
In our problem, we could call ELPA or ScaLAPACK for the first several SCF iterations without
estimating the DOS.
Since the intermediate band matrix has narrow semibandwidth, $n_{bw}=64$,
the cost of $LDL^H$ factorizations is small and
the spectrum can also be partitioned by using Sylvester's law of inertia.
To reduce storage cost, we store the banded matrix in LAPACK compact form.

After the first several SCF iterations, we have prior knowledge of eigenvalues, and
apply a \emph{k}-means clustering method to the previous eigenvalues and
partition the spectrum, just as done in~\cite{kecceli2018siesta,williams2020shift}.
Since the eigenvalues are stored as a one dimensional increasing array, the data structure is very simple, we propose
a simplified \emph{k}-means algorithm in the following subsection, which turns out to be good for the load balance problem.
In this work, the aim of using \emph{k}-means is to find a good partition of the spectrum,
i.e., separate the clustered eigenvalues.

\subsubsection{\emph{K}-means Clustering}

\begin{figure}[t]
\centering
\includegraphics[width=4.3in,height=2.2in]{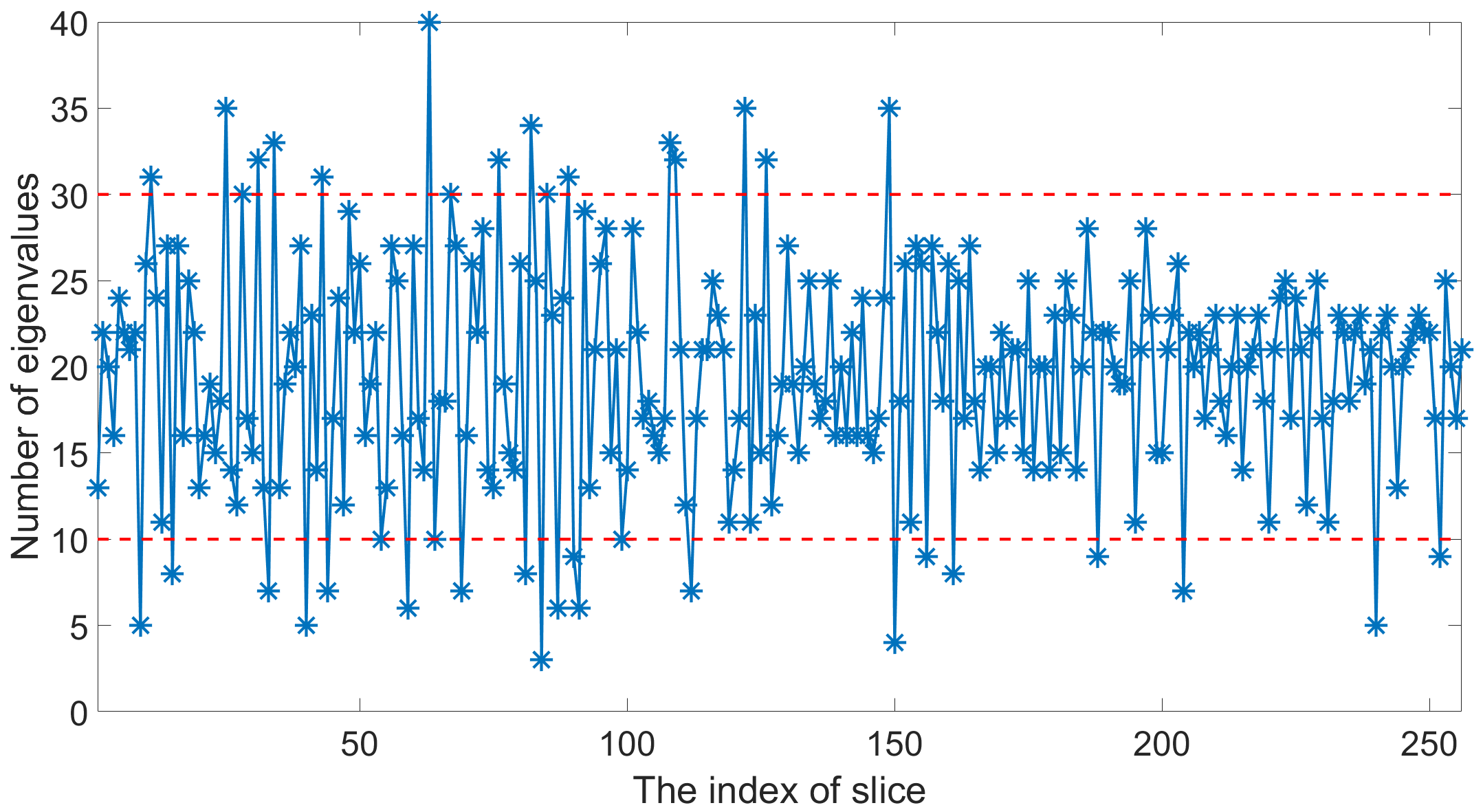}
\Description{}
\caption{The number of eigenvalues of 256 slices for \texttt{H2O}.}
\label{fig:imbalance}
\end{figure}

When we have a prior knowledge of the distribution of eigenvalues,
we can use one clustering algorithm to divide the spectrum into different slices
based on the previous computed eigenvalues.
One famous clustering algorithm is the \emph{k}-means algorithm~\cite{lloyd1982least}.
The goal is to group any cluster of eigenvalues
into a separate slice and let the boundaries of slices not cross any cluster of eigenvalues.
Figure~\ref{fig:kmeans} shows the boundaries of slices when partitioning $256$ eigenvalues
into eight groups, where the red vertical lines mark the boundaries of slices.
From it we can see that the \emph{k}-means algorithm can capture the clusters of eigenvalues and
separate them well.
One potential problem of \emph{k}-means partitioning is that different slices may be quite imbalanced,
i.e., some slices may have many eigenvalues while others may have very few eigenvalues,
just as shown in Figure~\ref{fig:imbalance}, especially when partitioning into many
slices.
Figure~\ref{fig:imbalance} shows the result of partitioning the first $5000$ eigenvalues
of matrix \texttt{H2O} from the SuiteSparse collection~\cite{davis2011university} into $256$ slices,
where the \emph{y}-axis represents the number of eigenvalues in each slice.
One slice contains only two eigenvalues while another one has 40 eigenvalues.
Finding the optimal clustering of eigenvalues is a difficult task.
In this work, we focus more on the goodness of separation
instead of the load imbalance problem.
The results of Example 2 in section~\ref{sec:num} show that the convergence rate of spectrum slicing
algorithms is more sensitive to the relative gap between eigenvalues in different slices
than the number of eigenvectors to be computed.

A good clustering should satisfy the following objectives, see also~\cite{kecceli2018siesta}:
\begin{enumerate}
  \item The clustered eigenvalues (separated by less than $10^{-7}$)
        should be kept in the same slice;

  \item The number of eigenvalues in different slices should be nearly equal;

  \item The eigenvalue range with each slice should be minimized.
\end{enumerate}
Objective $1$ is important for the orthogonality of computed eigenvectors,
and therefore should be considered first.
Objectives $2$ and $3$ concern the load imbalance problem, which is considered
secondary in this paper.

The classical \emph{k}-means algorithm is designed for multi-dimensional and
complicated data sets, and does not exploit the simplicity in the case of a one dimensional point set.
Our problem is very simple.
There is a one dimensional set of points $\mathcal{C}=\{d_1,d_2,\cdots,d_n\}$,
which are ordered increasingly, and we
want to partition $\mathcal{C}$ into $K$ groups in the hope that each group has nearly equal
number of points.
We propose a simplified and fast one dimensional \emph{k}-means algorithm,
which has the following advantages comparing with traditional \emph{k}-means algorithm.
\begin{enumerate}
  \item This method has less computational complexity, which only requires $n+K$ floating pointing operations (flops),
  and $n$ is the number of points and $K$ is the number of groups. Note that classical \emph{k}-means algorithm requires $(K+1)n$ flops.


  \item This method usually gets a partition with smaller variance of numbers of points in all slices.
\end{enumerate}
Our algorithm solves the same problem as the classical \emph{k}-means method,
\begin{equation}
\label{eq:kmeans_obj}
J=\min_{\mu_j} \sum_{i=1}^n \sum_{j=1}^K \|d_i-\mu_j\|^2,
\end{equation}
where $\mu_j, j=1,2,\cdots,K$ are the centroid of each group
$G_1, G_2, \cdots, G_K.$
The whole process is shown in Algorithm~\ref{alg:kmeans}.
In our implementation, we use Euclidean norm and
the centroid of group $G_i$ is computed as
$\mu_i=(d_l+d_r)/2$, where $d_l$ and $d_r$ are the
leftmost and rightmost points in $G_i$, respectively.
Algorithm~\ref{alg:kmeans} is allowed to adaptively adjust the parameter $K$ during iterations.
Therefore, the return number of slices may be less than the input value
$K$. When it happens, we will randomly generate $K$ centroids and
run the classical \emph{k}-means algorithm to generate $K$ partitions when
the number of slices is given.
During all our experiments, Algorithm~\ref{alg:kmeans} always returns the same
number of slices as given by $K$.

The spectrum is partitioned by recording the number $m_i$ of each slice,
and the boundary set $\mathcal{B}$ is updated in the next SCF iteration,
$b_i = \frac{d_{j_1}+d_{j_2}}{2}$, where $d_{j_1}$ is the largest value
in $G_{i}$ and $d_{j_2}$ is the smallest value in $G_{i+1}$.
Since Algorithm~\ref{alg:kmeans} only costs linear time, we perform
\emph{k}-means clustering at every SCF iteration.
But, we find it is not necessary.
For the last few iteration steps of SCF before convergence,
the spectrum partitions computed by Algorithm~\ref{alg:kmeans} are the same.
When the spectrum is repartitioned, the number of eigenvalues in each slice
may become different. It happens in the first few SCF iterations.
When $m_i$ is larger than before, the previous initial vectors $\widehat{X}_{m_i}$ can be interlaced
with some random vectors.
When $m_i$ is smaller than before, some previous initial vectors can be removed.

We compare Algorithm~\ref{alg:kmeans} with two other partition approaches in section~\ref{sec:num}.
One is dividing the eigenvalues equally into $K$ groups and every group
has nearly the same number of eigenvalues, $nev/K$, named by \texttt{Approach-1},
and the boundaries of slices are computed as Algorithm~\ref{alg:kmeans},
$b_i = \frac{d_{j_1}+d_{j_2}}{2}$, where $d_{j_1}$ is the largest value
in $G_{i}$ and $d_{j_2}$ is the smallest value in $G_{i+1}$.
\texttt{Approach-1} has perfect load balance since every slice has nearly
the same number of eigenvalues when used in \texttt{PDESHEP}.
However, the boundary $b_i$ may be inside a cluster of eigenvalues, i.e.,
there are many eigenvalues near the boundary $b_i$, and it may
introduce orthogonality problem, and slow down the convergence.
The other one, \texttt{Approach-2}, is dividing the spectrum interval into $K$ pieces equally,
which is rough.
\texttt{Approach-2} may have both load imbalance and orthogonality problems, and is
not suggested to be used in practice.

\SetNlSty{texttt}{(}{)}
\begin{algorithm}[thb]
\SetAlgoLined
\KwIn{Dataset $\mathcal{C}=\{d_1,d_2,\cdots,d_n\}$ in increasing order, and the number $K$ of clusters.}
\KwOut{The boundary point set $\mathcal{B}=\{b_0,b_1,\cdots,b_{K}\}$, and
the number set of points $\mathcal{M}=\{m_1,m_2,\cdots,m_K\}$.}

\nl Partition $\mathcal{C}$ into $K$ groups $\{G_i\}_{i=1}^K$ such that
$\left( n-\lfloor n/K \rfloor \cdotp K\right)$ groups have $\left( \lfloor n/K \rfloor+1 \right)$ points and the other groups have $\lfloor n/K \rfloor$ points\;

\nl Construct $\mathcal{B}=\{b_i\}_{i=0}^K$, where $b_i = \frac{d_{j_1}+d_{j_2}}{2}$, where $d_{j_1}$ is the largest in $G_{i}$ and $d_{j_2}$ is the smallest in $G_{i+1}$ for $i=1,\cdots,K-1$, $b_0$ and $b_{K}$ are the lower and upper bounds of $\mathcal{C}$, respectively\;

\nl Let $\mathcal{B}_0=\{\empty \}$\;

\nl \While{$(\mathcal{B} \ne \mathcal{B}_0)$
    }{
      $\mathcal{B}_0 =\mathcal{B}$\;
      \For{ $(j=1; j \le K-1; j++)$
       }{
\nl         $\ell = \textbf{argmin}_{d_\ell \in \mathcal{C}} \left\{ d_\ell > \text{centroid}(\text{centroid}(G_i), \text{centroid}(G_{i+1})) \right\},$
         where \emph{centroid}$(\cdotp)$ is the operator of calculating the centroid\;
\nl        $a_i = (d_{\ell}+d_{\ell-1})/2$\;
\nl         update the boundary between $G_i$ and $G_{i+1}$, and record
         the number $m_i$ of points in $G_i$\;
       }
       \While{ $(\exists b_p, b_q \in \mathcal{B},$ s.t. $b_p == b_q)$
       }{
\nl        $\mathcal{B} = \mathcal{B}\backslash b_q$, $K=K-1$\;
       }
        }
\Return{ $\mathcal{B}=\{b_i\}_{i=0}^K, \mathcal{M}=\{m_i\}_{i=1}^K$.}

\caption{A simplified determinate one dimensional \emph{k}-means algorithm}
\label{alg:kmeans}
\end{algorithm}

\subsection{Validation and Other Details}

The contour integral based algorithm can compute all the wanted eigenvalues enclosed in a given interval.
To guarantee that $nev$ eigenpairs can be computed, the lower bound $b_0$ and upper bound $b_K$
must enclose all the wanted $nev$ eigenvalues\footnote{For simplicity, we are assuming that
the smallest \emph{nev} eigenvalues are required. The cases for the largest or interior \emph{nev}
eigenvalues are similar. The users should provide how many (\emph{nev}) eigenvalues are required.}.
First, the upper bound $b_K$ should be larger than the $nev$-th eigenvalue
of the current matrix.
Since we know the previous eigenvalues, $b_K$ can be computed by $b_K=(1+\alpha)d_{nev}$,
where $\alpha$ is a small real number and $d_{nev}$ is the
 $nev$-th previous eigenvalue.
Note that $b_K$ should not be too large. Otherwise, $[b_0, b_K]$ will enclose too many eigenvalues.
$b_K$ can be checked by computing the inertia of $(D^{(\ell)}-b_k I)$.
If $b_K$ is not large enough, we iteratively enlarge it by ratio $\alpha$, $b_K = (1+\alpha) b_K$,
and check it again until there are more than $nev$ eigenvalues less than $b_K$.
Algorithm~\ref{alg:lower_upper} illustrates process of computing the upper boundary $b_K$.
Since LAPACK does not support $LDL^H$ factorization for banded symmetric matrices,
we implemented a simplified version of the retraction algorithm~\cite{kaufman2007retraction}, which
is used to only compute the inertia without recording the matrix $L$, and
the storage and computational cost are further reduced, becoming very small.

Another approach is to enlarge $b_K$ as $b_K = b_K + \hat{m}_1 \delta$, where
$\hat{m}_1$ is an integer and $\delta = d_{nev} - d_{nev-1}$, the distance between the
$nev$-th and $(nev-1)$-th previous eigenvalues.
When $b_k$ is not large enough, it will be enlarged by another $\hat{m}_1\delta$.
Its computation is similar to the first approach.
Since we are solving sequences of eigenvalue problems, the parameters $\alpha$ and $\hat{m}_1$
can be updated during the SCF iterations.
In our implementation, we let $\alpha=\pm 10^{-3}$ and $\hat{m}_1=6$.
For the first fewer SCF iterations, Algorithm~\ref{alg:lower_upper} needs to iterate at most $2$--$3$ times
to obtain the upper boundary point.
Since $D^{(\ell)}$ is banded with small bandwidth, the computation of these inertia is very cheap.

The lower bound $b_0$ can be computed similarly, $b_0=(1+\beta)d_1$, and make sure
$b_0 < d_1$ is a lower bound.
When computing the smallest few eigenvalues, the parameter $\beta$ is relatively easy to choose
since any lower bound works.
Similarly, $b_0$ can also be updated by $b_0 = b_0 + \hat{m}_2 \delta_2$, where
$\hat{m}_2$ is integer and $\delta_2=d_1 -d_2$, negative.
When computing interior eigenvalues, $\beta$ should be chosen carefully in case
$b_0$ is too small to enclose too many eigenvalues in $[b_0, b_K]$.
Since the process is similar to Algorithm~\ref{alg:lower_upper}, it is not included.

%
%
%
%
%

\SetNlSty{texttt}{(}{)}
\begin{algorithm}[thb]
\SetAlgoLined
\KwIn{The  $nev$-th previous eigenvalue $d_{nev}$ and integer $nev$, parameters $\alpha$, $\hat{m}_1$ and $\delta$
     which equals to $d_{nev}-d_{nev-1}$, and matrix $D^{(\ell)}$.}
\KwOut{The upper boundary point ${b}_K$.}

\nl Compute ${b}_K = d_{nev}+\max( \alpha \cdotp d_{nev}, \hat{m}_1 \delta)$\;

\nl Let $flag$=true and $it=1$\;

 \While{$(flag)$
    }{
\nl      Compute the inertia of $(D^{(\ell)}-{b}_K I)$, $(n_{-}, n_0, n_+)$\;
      \uIf{ $(n_{-} \le nev)$
       }{
\nl        $flag$=false\;
       }
       \uElse{
\nl         ${b}_K = {b}_K + \max( \alpha \cdotp d_{nev}, \hat{m}_1 \delta)$\;
\nl         $it = it +1$\;
       }

        }
\Return{ ${b}_K$.}

\caption{Compute the upper boundary point $b_K$.}
\label{alg:lower_upper}
\end{algorithm}

The spectrum is divided into $K$ slices and one process group works on one
slice.  Note that only the first process group needs to compute the lower boundary point $b_0$,
and only the $K$-th process group needs to compute the upper boundary point $b_K$.
The $i$-th process group computes the eigenvalues in the interval $(b_{i-1}, b_i)$, and
every process can compute the exact number $m_i$ of eigenvalues in the current interval by
computing the inertia of matrices $(D^{(\ell)}-b_{i-1} I)$ and $(D^{(\ell)}-b_i I)$.
Therefore, every process group can check whether the computed eigenvalues are
correct or not, and if some eigenvalues are missed, the corresponding process group should rerun
FEAST by using more starting vectors.

In our problem, the value $m_i$ is the exact number of eigenvalues in the $i$-th slice, which
is used to check the correctness of computed eigenpairs.
For stability, we further enlarge the number of starting vectors in FEAST,
$m_i^{(0)}=\max(1.3*m_i,m_i+10)$\footnote{The number $10$ is chosen arbitrarily in order to deal with
the case that $m_i$ is too small and $m_i^{(0)}$ equals $1.3 \times m_i$.}.
In our experiments, FEAST always computes more than $nev$ eigenvalues and no
eigenvalues are missed. 

%
%
%
%

\section{Numerical results}
\label{sec:num}

All the experimental results are obtained on the Tianhe-2 supercomputer~\cite{Liao-TH2, Liao-HPCG},
located at Guangzhou, China. Each compute node is equipped with two $12$-core Intel Xeon E5-2692 v2 CPUs.
The details of the test platform and environment of the compute nodes are shown in Table~\ref{tab:tianhe2a}.
For all these numerical experiments, we only used plain MPI, running $24$ MPI processes per node in principle, that is,
one process per core.
Our codes are written in Fortran 90 and part of our codes will
be released on Github\footnote{\url{https://github.com/shengguolsg/}}.

\begin{table}[ptbh]
\caption{The test platform and environment of one node.}
\label{tab:tianhe2a}
\begin{center}%
\begin{tabular}
{|c|c|}\hline
Items  & Values  \\ \hline
2*CPU  & Intel Xeon CPU E5-2692 v2@2.2GHz \\ \hline
Memory size & 64GB (DDR3) \\ \hline
Operating System  & Linux 3.10.0 \\ \hline
Compiler &Intel ifort \\ \hline
Optimization & -O3 -mavx
\\ \hline
\end{tabular}
\end{center}
\end{table}

\textbf{Example 1.}
For the experiments in this example, we use an eigenvalue problem obtained from a
DFT simulation using the FLEUR code~\cite{fleur-web-page} that produces a matrix labeled \texttt{NaCl}, which
has been used in~\cite{winkelmann2019chase}.
This matrix is Hermitian and its size is $9273$.
We use ELPA to compute partial and full eigendecomposition of this matrix,
and the point is to show that step (II) has poor scalability.
The codes of ELPA that we used are obtained from package ELSI v2.5.0~\cite{yu2020elsi}.
The performance of ELPA highly depends on the parameter $n_b$, the block size used
for storing matrices in 2D BCDD format.
Many works suggest that $n_b=128, 64, 32$ are usually good choices, and smaller $n_b$ such as
less than $16$ will degrade the performance, see~\cite{Elpa} and~\cite{thpc_pbsdc} for more details.
There are no big differences for ELPA when choosing $n_b=128, 64$ or $32$ on Tianhe-2
supercomputer.
All the performance results of ELPA are obtained by choosing $n_b=64$ in this paper.

The results for computing $nev=2000, 5000$ and $9273$ eigenpairs
are shown in Table~\ref{tab:elpa_step4} when using
$64, 256, 1024$ and $4096$ processes.
We use matrix \texttt{NaCl} to construct a standard eigenvalue problem, solve it
by ELPA, and measure
the times of step (I)-(III). From the results, we can see that
\begin{itemize}
  \item the timings of step (II) nearly stop decreasing after $n_p \ge 1024$;
  \item step (II) takes about half of the total time, especially when $n_p$ is large.
\end{itemize}
Therefore, it is better to replace step (II) by a more efficient and faster solver.
In this paper, we use a spectrum slicing approach for this step.
The next example shows that this approach can reduce the times of step (II) by half, and
our final algorithm (\texttt{PDESHEP}) can be more than $1.25$ times faster than ELPA.

\begin{table}
\small
\caption{The execution times (in second) of each step of ELPA when choosing different $nev$.}
\label{tab:elpa_step4}
\begin{center}
\begin{tabular}[c]{|c|c|c|c|c|c|c|c|c|c|} \hline
 \multirow{2}{*}{$n_p$} & \multicolumn{3}{c}{$nev=2000$} & \multicolumn{3}{|c|}{$nev=5000$} & \multicolumn{3}{|c|}{$nev=9273$} \\ \cline{2-10}
  &  Step (I) & Step (II) & Step (III) & Step (I) & Step (II) & Step (III) & Step (I) & Step (II) & Step (III)\\ \hline
64   & $7.62$ & $4.49$ & $1.82$ &  $7.67$ & $8.80$ & $4.49$  & $7.57$ & $15.10$ & $8.03$\\ \hline
256   &  $3.47$ & $1.96$ & $0.68$  & $3.44$ & $3.01$ & $1.50$  & $3.45$ & $4.78$ & $2.70$ \\ \hline
1024  &  $2.09$ & $2.05$ & $0.33$  & $2.17$ & $2.19$ & $0.62$ & $2.13$ & $2.83$ & $1.07$  \\ \hline
4096 &  $1.85$ & $2.18$ & $0.32$  & $1.93$ & $2.03$ & $0.50$  & $1.92$ & $2.30$ & $0.54$ \\ \hline
\end{tabular}
\end{center}
\end{table}

\begin{figure}[b]
\centering
\includegraphics[width=4.5in,height=1.5in]{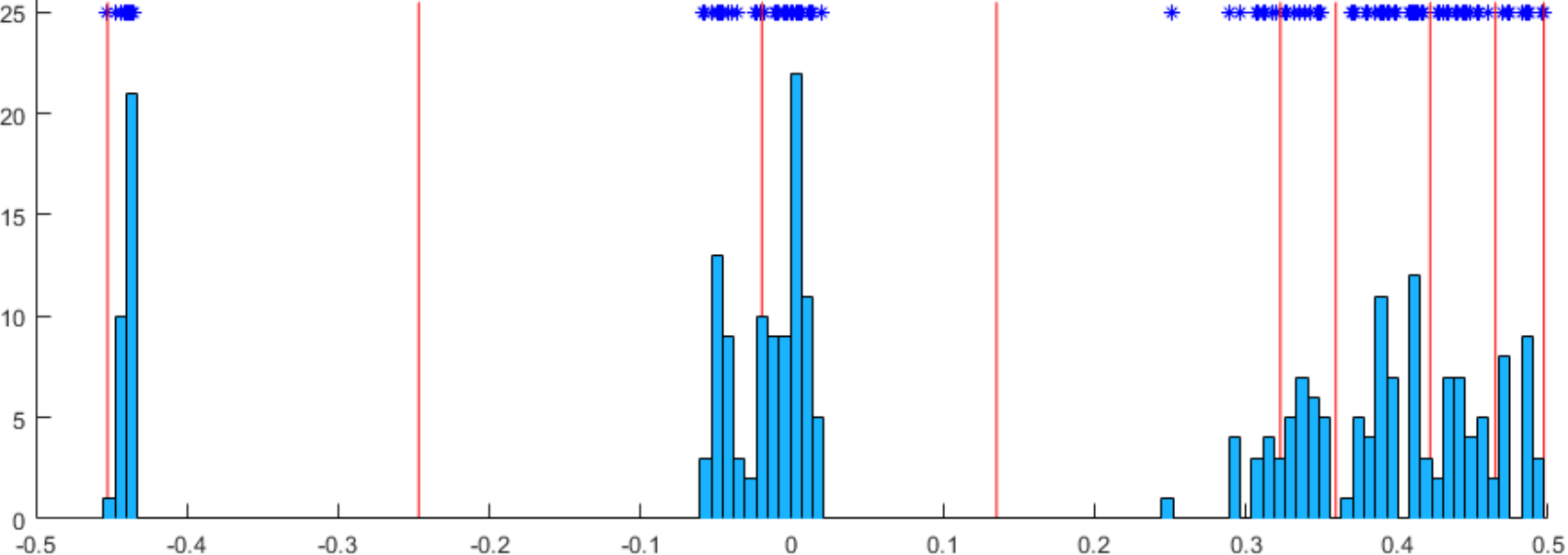}
\Description{The procedure of our algorithm}
\caption{Using \emph{k}-means to clustering 256 eigenvalues of matrix \texttt{NaCl} into 8 groups. The red lines are the boundaries of slices, and
the histogram shows the number of eigenvalues in each slice. }
\label{fig:kmeans}
\end{figure}

\textbf{Example 2.} In this example, we use a sequence of eigenvalue problems.
It is the \texttt{NaCl} sequence of size $9273$, used in the previous example and in~\cite{winkelmann2019chase}.
There are $16$ such matrices.
We first use \texttt{PDESHEP} to compute the partial eigendecomposition of the $15$th matrix
by using the eigenpairs of the $14$th and the purpose is to show that
\texttt{PDESHEP} works for highly related eigenvalue problems.
The first $256$ eigenvalues of the $14$th matrix are shown in Figure~\ref{fig:kmeans}.
Since the SCF loop is nearly converged, eigenvalues of the $15$th and $14$th matrices are quite similar and
the starting vectors obtained from the $14$-th matrix are also very good approximations.

The ChASE algorithm~\cite{winkelmann2019chase} computes $nev=256$ smallest eigenpairs of \texttt{NaCl}
by using Chebyshev accelerated subspace iterations.
When computing many eigenpairs, ChASE would become much slower since it does not use
spectrum slicing strategy.
The advantage of our algorithm is that it can compute many (even all) eigenpairs efficiently.
In this example, we compute $nev=2000, 5000$ and $9273$ eigenpairs for this
matrix.

\begin{figure}[ptbh]
\centering
\subcaptionbox{The speedup over ELPA for Stage II.\label{fig:nacl_spd_p3}}{
\includegraphics[width=2.5in,height=2.3in]{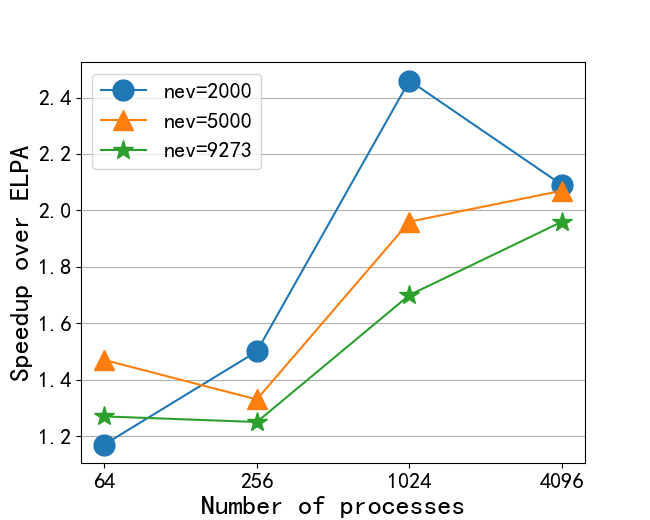}
}
\subcaptionbox{The total speedup of \texttt{PDESHEP} over ELPA.\label{fig:nacl_totalspeedup}}{
\includegraphics[width=2.5in,height=2.3in]{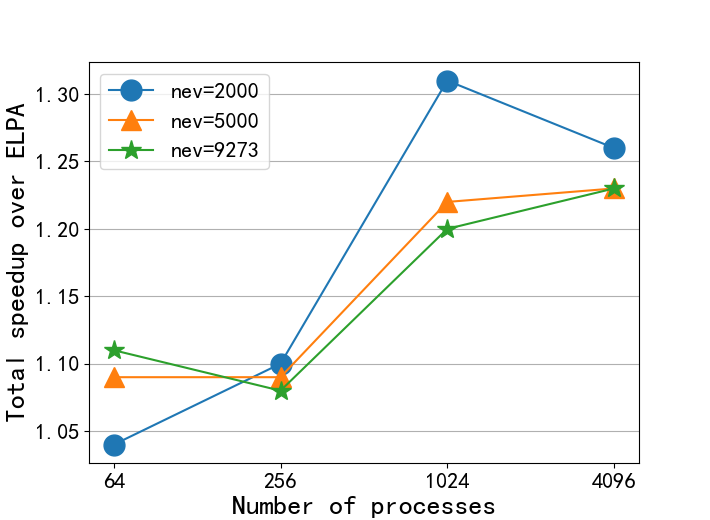}
}
\caption{The speedup of \texttt{PDESHEP} over ELPA.}
\label{fig:Nacl_speed}
\end{figure}

Let the block size $n_b$ be $64$, which is nearly optimal on our test platform.
We first reduce the \texttt{NaCl} matrix to a banded Hermitian matrix
with semi-bandwidth $n_{bw}=64$, and call FEAST to compute the eigenpairs in all slices simultaneously.
Since step (I) and step (III) are the same as in ELPA, we mainly compare step (II) of ELPA and \texttt{PDESHEP}.
The timing results of FEAST shown in Fig.~\ref{fig:Nacl_speed} are obtained when using the
eigenvectors of the banded matrix from the previous SCF iteration.
If using random starting vectors, FEAST would require more time as illustrated
in Table~\ref{tab:wo_vectors}.
From the results in Table~\ref{tab:wo_vectors}, we can see that using the previous eigenvectors as starting vectors can be more than four times
faster than using random starting vectors when using $4096$ processes.
The results show that using previous eigenvectors as initial vectors provide a significant benefit when
the eigenvalue problems are highly correlated.
\texttt{PDESHEP} is suggested to be used for the last few SCF
iterations before convergence.

\begin{table}
\caption{The times (in second) of FEAST with and without initial vectors when applying on the $15$-th matrix.}
\label{tab:wo_vectors}
\begin{center}
\begin{tabular}[c]{|c|c|c|c|c|c|c|} \hline
 \multirow{2}{*}{$n_p$} & \multicolumn{2}{c}{$nev=2000$} & \multicolumn{2}{|c|}{$nev=5000$} & \multicolumn{2}{|c|}{$nev=9273$} \\ \cline{2-7}
  &  With & Without & With & Without & With & Without \\ \hline
64   & $3.75$ & $5.43$ & $5.86$  & $8.74$ & $ 11.70$ & $18.51$ \\ \hline
256   & $1.17$ & $1.96$ & $2.10$ & $3.27$ & $3.64$ & $5.15$ \\ \hline
1024  & $0.60$ & $0.86$  & $1.61$ & $1.99$ & $1.42$ & $2.24$   \\ \hline
4096 &  $0.87$ & $3.14$ & $0.62$ & $2.90$ & $0.71$ & $2.86$ \\ \hline
\end{tabular}
\end{center}
\end{table}

When comparing with step (II) of ELPA, we add the time of FEAST with the time
spent in data redistributions.
The speedups over Stage II of ELPA are shown in Fig.~\ref{fig:nacl_spd_p3}
when using different processes ($64, 256, 1024, 4096$) to compute $2000, 5000$ or $9273$ eigenpairs.
The speedup can be up to $2\times$ when using $4096$ processes.
The total speedups of \texttt{PDESHEP} over ELPA are shown in Fig.~\ref{fig:nacl_totalspeedup}.
\texttt{PDESHEP} can be more than 20\% faster than ELPA
when computing full or partial eigendecomposition and using $4096$ processes.

During the first few SCF iterations, the matrices can be quite different and
so do the eigenvalues and eigenvectors.
We now apply \texttt{PDESHEP} on a matrix on the early stages of SCF iterations and compare \texttt{PDESHEP}
with ELPA.
We choose the $10$-th \texttt{NaCl} matrix and use the same parameter $\epsilon=10^{-11}$.
We first compute the eigendecomposition of the $10$-th \texttt{NaCl} matrix and
using its eigenvalues to partition the spectrum of the $11$-th matrix.
Then, we use eigenvectors of the $10$-th matrix as initial vectors of FEAST.
We see that the `good` estimated initial vectors do help the convergence of FEAST but they
are not as significant as that obtained on the $15$-th matrix.
The results are shown in Table~\ref{tab:firststage}, where
the column \emph{With} shows the times (in second) cost by step (II) of FEAST when the initial vectors are constructed from
the eigenvectors of the $5$-th matrix;
the column \emph{Without} shows the times when using random vectors as initial vectors;
the column \emph{ELPA} shows the times cost by step (II) of ELPA.
For this matrix, we use at most $1024$ processes.
From the results, we can see that \texttt{PDESHEP} can still be
faster than ELPA for the $11$-th matrix.

\begin{table}
\caption{The times (in second) of FEAST with and without initial vectors when applying on the $11$-th matrix with $\epsilon=10^{-11}$.}
\label{tab:firststage}
\begin{center}
\begin{tabular}[c]{|c|c|c|c|c|c|c|c|c|c|} \hline
 \multirow{2}{*}{$n_p$} & \multicolumn{3}{c}{$nev=2000$} & \multicolumn{3}{|c|}{$nev=5000$} & \multicolumn{3}{|c|}{$nev=9273$} \\ \cline{2-10}
  &  With & Without & ELPA & With & Without & ELPA & With & Without & ELPA \\ \hline
64   & $3.22$ & $4.12$  & $4.49$ & $6.91$  & $10.25$ & $8.80$ & $14.58$ & $18.97$ & $15.10$  \\ \hline
256   & $1.57$ & $2.16$ & $1.96$ & $2.80$ & $3.72$   & $3.01$ & $4.66$  & $5.75$  & $4.78$ \\ \hline
1024  & $1.19$ & $1.84$ & $2.05$ & $1.74$ & $2.24$   & $2.19$ & $2.56$  & $3.45$  & $2.83$  \\ \hline
\end{tabular}
\end{center}
\end{table}

For the matrices appeared at the even earlier stages of SCF, for example the $5$-th matrix,
the 'good' estimated initial vectors and eigenvalues help very little.
FEAST may cost more than step (II) of ELPA when computing the eigenpairs
accurately with $\epsilon=10^{-11}$.
However, at the early stages of SCF iterations, the eigenvalue problems only need to be solved
approximately.
Therefore, we let $\epsilon=10^{-9}$ and find that
step (II) of \texttt{PDESHEP} is also faster than that of ELPA,
and the results are shown in Table~\ref{tab:firststage2} for the $6$-th matrix with $\epsilon=10^{-9}$.

\begin{table}
\caption{The times (in second) of step (II) for FEAST with and without initial vectors when applying on the $6$-th matrix with $\epsilon=10^{-9}$.}
\label{tab:firststage2}
\begin{center}
\begin{tabular}[c]{|c|c|c|c|c|c|c|c|c|c|} \hline
         \multirow{2}{*}{$n_p$} & \multicolumn{3}{c}{$nev=2000$} & \multicolumn{3}{|c|}{$nev=5000$} & \multicolumn{3}{|c|}{$nev=9273$} \\ \cline{2-10}
          &  With & Without & ELPA & With & Without & ELPA & With & Without & ELPA \\ \hline
        64   & $3.35$ & $4.38$  & $4.49$ & $7.22$  & $10.18$ & $8.80$ & $13.12$ & $18.66$ & $15.10$ \\ \hline
        256   & $1.30$ & $1.76$ & $1.96$ & $2.84$ & $3.44$   & $3.01$ & $4.16$ & $5.73$ & $4.78$ \\ \hline
        1024  & $1.18$ & $1.91$ & $2.05$  & $1.88$ & $2.84$  & $2.19$ & $2.46$ & $3.57$ & $2.83$  \\ \hline
        \end{tabular}
\end{center}
\end{table}

\begin{figure}[h]
\centering
\includegraphics[width=3.3in,height=2.5in]{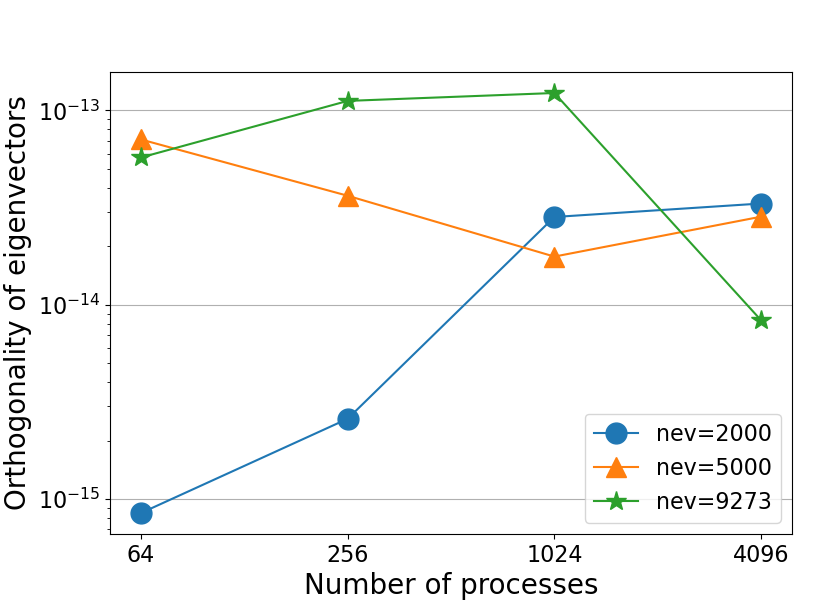}
\Description{}
\caption{The orthogonality of computed eigenvectors of matrix \texttt{NaCl}.}
\label{fig:nacl_orth}
\end{figure}

In this paper, the orthogonality is measured as
\begin{equation}
\mathrm{Orth} =\frac{\|X^H X-I\|_{\max}}{n},
\end{equation}
where $n$ is the dimension of the matrix, $X$ is the computed eigenvector matrix and
$I$ is the identity matrix of size $nev$.
The orthogonality of computed eigenvectors is shown in Fig.~\ref{fig:nacl_orth},
being about or less than $10^{-13}$ in all cases, an acceptable value for most material simulations.
Though the orthogonality of eigenvectors computed by \texttt{PDESHEP} is
not as accurate as those obtained by direct methods, which is usually
around the machine precision, it should be accurate enough for the SCF iterations,
The computed residuals of eigenpairs are all less than $10^{-11}$ since the
residual tolerance of FEAST is chosen to be $\epsilon =10^{-11}$ in this paper.
In the previous work~\cite{aktulga2014parallel}, the residual tolerance was $10^{-10}$, a little
larger than that used in this paper.

Since ChASE does not implement the spectrum slicing techniques, it is suggested to compute less than
$10$\% eigenpairs of matrix~\cite{aktulga2014parallel}.
We use ChASE to compute the first $1000$ eigenpairs of the $15$-th \texttt{NaCl} matrix, and
the orthogonality of computed eigenvector by ChASE is around $10^{-18}$ which is more accurate than
those computed by \texttt{PDESHEP}, but ChASE costs much more time than \texttt{PDESHEP}.
If we let $\epsilon = 10^{-12}$, the orthogonality of computed eigenvectors by \texttt{PDESHEP}
are also improved, around $10^{-16}$, when using $64$, $256$ and $512$ MPI processes.
The results are shown in Table~\ref{tab:orthogonality}.
The convergence of FEAST depends on the parameter $\epsilon$, and we find
it may take FEAST a lot of iterations to converge  if
$\epsilon < 10^{-13}$.

\begin{table}
\caption{The orthogonality of eigenvectors computed by ChASE and \texttt{PDESHEP} for the $15$-th \texttt{NaCl} matrix.}
\label{tab:orthogonality}
\begin{center}
\begin{tabular}[c]{c|ccc} \hline
 Method & {$n_p=64$} & {$n_p=256$} & {$n_p=512$} \\ \hline
\texttt{PDESHEP}   & $9.49e$-$16$ & $9.48e$-$16$ & $9.40e$-$16$ \\ \hline
ChASE   & $2.22e$-$18$ & $1.41e$-$18$ & $1.84e$-$18$ \\ \hline
\end{tabular}
\end{center}
\end{table}

The results of this example show that \texttt{PDESHEP} can compute both
partial and full eigendecomposition efficiently and the computed eigenvectors
are orthogonal with acceptable accuracy when the slices are constructed carefully.
From our experimental results, we claim that the \emph{k}-means based spectrum
partition strategy (Algorithm~\ref{alg:kmeans}) is helpful for the orthogonality of computed eigenvectors.
We compare Algorithm~\ref{alg:kmeans} with two other partition approaches:
\texttt{Approach-1} and \texttt{Approach-2},
and Algorithm~\ref{alg:kmeans} performs the best.

\textbf{Example 3.} In this example, we use some relatively large matrices to
test our algorithm.
The first class is \texttt{AuAg} with dimension $n=13,379$, which originates from simulations using FLAPW~\cite{jansen1984total}
and has been used in~\cite{winkelmann2019chase}. They have different spectral distribution compared to the \texttt{NaCl} dataset.
There are $25$ such matrices, and we compute the eigendecomposition of the $24$th matrix
by using the information from the $23$th matrix.

The other two matrices are from the SuiteSparse collection~\cite{davis2011university}.
One is \texttt{Si5H12} with $n=19,896$ and $nnz=379,247$, and the other is
the \texttt{H2O} matrix with $n=67,024$ and $nnz=2,216,736$,  which are both real symmetric and sparse.
Since there is only one such matrix for each type, we construct another two matrices by adding a small
random symmetric perturbation
to the nonzero entries of these matrices such that they have the same structure and similar eigenvalues.
The nonzero entry $a_{i,j}$ is perturbed as
\begin{equation}
a_{i,j} = a_{i,j}(1+\tau \eta_{i,j} ),
\end{equation}
where $\tau=10^{-4}$ and $\eta_{i,j}$ is a random real number in the interval $[0, 1]$.
Note that the contour integral based algorithms always operate with complex arithmetic.
Therefore, FEAST may require more flops than other iterative methods for real matrices.

In this example, we compute the first $5000$ smallest eigenpairs and use $256, 1024$ and $4096$ processes,
respectively.
The results of ELPA are shown in Table~\ref{tab:large}.
For these experiments, the convergence residual tolerance is also $\epsilon =10^{-11}$, and
the dimension of the subspace is $m_i^{(0)}=\max( 1.3\times m_i, m_i+10)$, where
$m_i$ is the estimated number of eigenvalues in the current slice.
We note that $m_i^{(0)}$ can affect the performance considerably.
When $m_i^{(0)}$ is small, it makes the linear equations with fewer right-hand sides much
easier to solve.
When $m_i^{(0)}$ is too small, it may make FEAST not converge.
There is a trade-off between speed and robustness.


\begin{table}
\small
\caption{The execution times (in second) of each step of ELPA for large matrices.}
\label{tab:large}
\begin{center}
\begin{tabular}[c]{|c|c|c|c|c|c|c|c|c|c|} \hline
 \multirow{2}{*}{$n_p$} & \multicolumn{3}{c}{\texttt{AuAg}} & \multicolumn{3}{|c|}{\texttt{Si5H12}} & \multicolumn{3}{|c|}{\texttt{H2O}} \\ \cline{2-10}
  &  Step (I) & Step (II) & Step (III) & Step (I) & Step (II) & Step (III) & Step (I) & Step (II) & Step (III) \\ \hline
256   &  $8.19$ & $5.47$ & $2.87$  & $6.84$ & $5.82$ & $2.08$  & $17.41$ & $67.03$ & $18.87$  \\ \hline
1024  &  $4.01$ & $3.56$ & $1.13$  & $4.13$ & $4.45$ & $1.21$  & $64.34$ & $24.72$ & $9.70$ \\ \hline
4096 &  $3.29$ & $4.17$ & $0.76$  & $3.79$ & $4.14$ & $1.28$   & $27.45$ & $17.53$ & $5.57$ \\ \hline
\end{tabular}
\end{center}
\end{table}

\begin{figure}[ptbh]
\centering
\subcaptionbox{The speedup over ELPA for Stage II.\label{fig:large_StepIV}}{
\includegraphics[width=2.6in,height=2.1in]{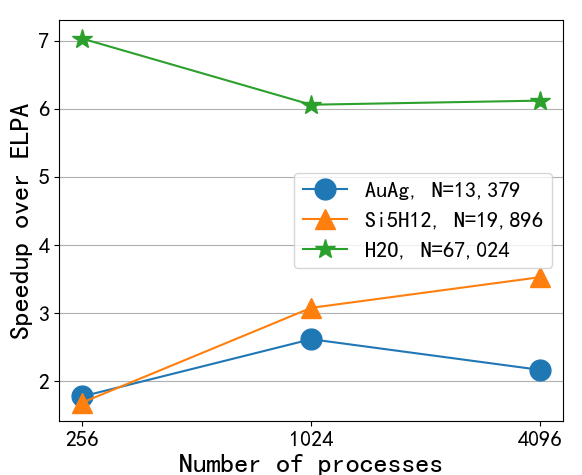}}
\subcaptionbox{The total speedup of \texttt{PDESHEP} over ELPA.\label{fig:large_totalspeedup}}{
\includegraphics[width=2.6in,height=2.05in]{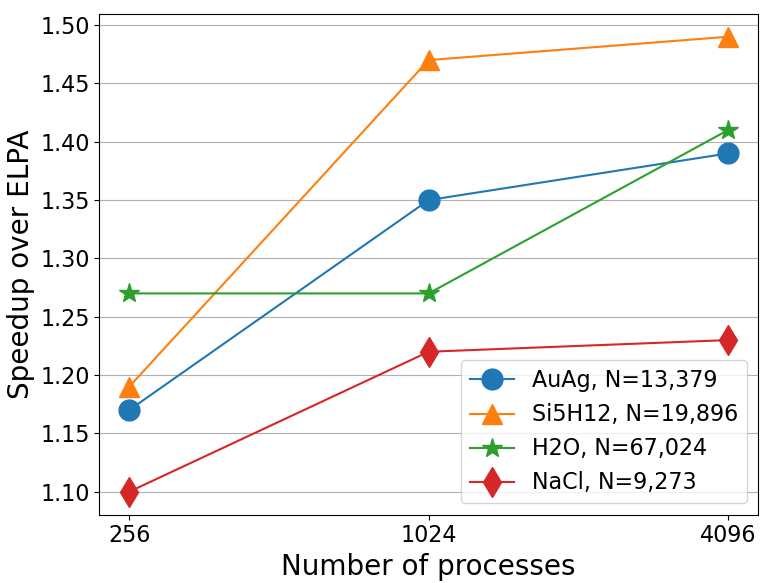}
}
\caption{The speedup of \texttt{PDESHEP} over ELPA for large matrices.}
\label{fig:Large_speed}
\end{figure}

The speedups of FEAST over step (II) of ELPA are shown in Fig.~\ref{fig:large_StepIV}.
The results in Fig~\ref{fig:large_StepIV} show that it is difficult to predict the speedup of
spectrum slicing algorithms over Stage II of ELPA, that is because the partition of the
spectrum may affect the performance of spectrum slicing algorithms.
It is also very hard to know how many slices to partition is best, which
depends on the distribution of eigenvalues and the used spectrum slicing algorithm.
Figure~\ref{fig:imbalance} shows the number of eigenvalues of the $256$ slices for \texttt{H2O} when using $4096$ processes.
It is not well balanced, and most slices have eigenvalues between $10$ and $30$, but
one slice only has $2$ eigenvalues and one slice
has $40$ eigenvalues.
Comparing with load imbalance problem, we find that good separation of eigenvalues may be more important.
If one boundary is inside a cluster of eigenvalues, it may introduce convergence and orthogonality problems.
We use our simplified \emph{k}-means algorithm to partition the spectrum,
and the load balance problem is also taken into account.
The total speedups of \texttt{PDESHEP} over ELPA are shown in Fig.~\ref{fig:large_totalspeedup},
which are always larger than $1.20$ and up to $1.49$ when using $4096$ processes.
The wall times of FEAST used to compute the second step, step (II), are shown
in Table~\ref{tab:jianxian}.

\begin{table}
\caption{The execution times of step (II) in second of \texttt{PDESEHP} for large matrices.}
\label{tab:jianxian}
\begin{center}
\begin{tabular}[c]{c|c|c|c} \hline
{$n_p$} & {\texttt{AuAg}} & {\texttt{Si5H12}} & {\texttt{H2O}} \\ \hline
256   & $3.08$ & $3.44$ & $9.53$  \\ 
1024  & $1.36$ & $1.44$ & $4.07$ \\ 
4096  & $1.92$ & $1.17$ & $2.86$ \\ \hline
\end{tabular}
\end{center}
\end{table}

\section{Conclusions and future works}
\label{sec:conclusion}

In this paper, we propose a new framework for computing partial or full eigendecomposition
of dense Hermitian or real symmetric matrices, which can be easily extended to generalized
eigenvalue problems.
This new framework is denoted by \texttt{PDESHEP}, which combines the classical direct methods with spectrum slicing algorithms.
In this work, we combine ELPA~\cite{Elpa} with FEAST~\cite{kestyn2016pfeast}.
\texttt{PDESHEP} can also be combined with other spectrum slicing algorithms such as those implemented in SLEPc~\cite{campos2012strategies} or EVSL~\cite{li2019eigenvalues}.
This will be a topic for future work.

In \texttt{PDESHEP}, a dense matrix is first reduced to a banded form and then spectrum slicing algorithms
are used to compute its partial eigendecomposition instead of further reducing it to tridiagonal form.
Therefore, the symmetric bandwidth reduction (SBR) process is avoided, which consists of memory-bounded
operations.
For self-consistent field (SCF) eigenvalue problems which are correlated,
\texttt{PDESHEP} can be about $1.25$ times faster than ELPA.
Numerical results are obtained on Tianhe-2 supercomputer and up to $4096$ processes are used.
The tested matrices include dense Hermitian matrices from real applications, and
large real symmetric sparse matrices downloaded from the SuiteSparse collection~\cite{davis2011university}.
Another future work is to incorporate \texttt{PDESHEP} into ELSI~\cite{yu2020elsi} or some electronic-structure calculation packages
such as QE~\cite{giannozzi2009quantum} to do some real simulations.

\begin{acks}
The authors would like to acknowledge many helpful discussions with Chao Yang from
LBNL in USA, and Chun Huang, Min Xie and Tao Tang from NUDT in China,  and
also want to thank Dr. Di Napoli from J\"{u}lich Supercomputing Centre for having provided some matrices for the eigenproblems used in the experiments.
Most of this work was done while the first author was visiting Universitat Polit\`{e}cnica de Val\`{e}ncia during 2020,
and he wants to thank J.~Roman for his hospitality.
This work is in part supported by National Key RD Program of China (2017YFB0202104), National Natural Science Foundation of China
(No. NNW2019ZT6-B20, NNW2019ZT6-B21, NNW2019ZT5-A10),
NSF of Hunan (No. 2019JJ40339), and
NSF of NUDT (No. ZK18-03-01).
Xinzhe Wu was in part supported by European Union (EU) under PRACE-6IP-WP8 Performance Portable Linear Algebra project (grant agreement ID 823767).
Jose E. Roman was supported by the Spanish Agencia Estatal de Investigaci\'{o}n (AEI) under project SLEPc-DA (PID2019-107379RB-I00).
\end{acks}

\bibliographystyle{ACM-Reference-Format}
\bibliography{../TPDS_PSDC/thesis-refer}

%

\end{document}